\documentclass[11pt]{article}
\usepackage{amssymb,amsmath,latexsym}
\usepackage{pdfsync}

\begin{document}

\begin{doublespace}

\def\1{{\bf 1}}
\def\ind{{\bf 1}}
\def\nn{\nonumber}

\def\sA {{\cal A}} \def\sB {{\cal B}} \def\sC {{\cal C}}
\def\sD {{\cal D}} \def\sE {{\cal E}} \def\sF {{\cal F}}
\def\sG {{\cal G}} \def\sH {{\cal H}} \def\sI {{\cal I}}
\def\sJ {{\cal J}} \def\sK {{\cal K}} \def\sL {{\cal L}}
\def\sM {{\cal M}} \def\sN {{\cal N}} \def\sO {{\cal O}}
\def\sP {{\cal P}} \def\sQ {{\cal Q}} \def\sR {{\cal R}}
\def\sS {{\cal S}} \def\sT {{\cal T}} \def\sU {{\cal U}}
\def\sV {{\cal V}} \def\sW {{\cal W}} \def\sX {{\cal X}}
\def\sY {{\cal Y}} \def\sZ {{\cal Z}}

\def\bA {{\mathbb A}} \def\bB {{\mathbb B}} \def\bC {{\mathbb C}}
\def\bD {{\mathbb D}} \def\bE {{\mathbb E}} \def\bF {{\mathbb F}}
\def\bG {{\mathbb G}} \def\bH {{\mathbb H}} \def\bI {{\mathbb I}}
\def\bJ {{\mathbb J}} \def\bK {{\mathbb K}} \def\bL {{\mathbb L}}
\def\bM {{\mathbb M}} \def\bN {{\mathbb N}} \def\bO {{\mathbb O}}
\def\bP {{\mathbb P}} \def\bQ {{\mathbb Q}} \def\bR {{\mathbb R}}
\def\bS {{\mathbb S}} \def\bT {{\mathbb T}} \def\bU {{\mathbb U}}
\def\bV {{\mathbb V}} \def\bW {{\mathbb W}} \def\bX {{\mathbb X}}
\def\bY {{\mathbb Y}} \def\bZ {{\mathbb Z}}
\def\R {{\mathbb R}} \def\RR {{\mathbb R}} \def\H {{\mathbb H}}
\def\n{{\bf n}} \def\Z {{\mathbb Z}}

\newcommand{\expr}[1]{\left( #1 \right)}
\newcommand{\cl}[1]{\overline{#1}}
\newtheorem{theorem}{Theorem}[section]
\newtheorem{lemma}[theorem]{Lemma}
\newtheorem{definition}[theorem]{Definition}
\newtheorem{proposition}[theorem]{Proposition}
\newtheorem{corollary}[theorem]{Corollary}
\newtheorem{remark}[theorem]{Remark}
\newtheorem{example}[theorem]{Example}
\numberwithin{equation}{section}
\def\ee{\varepsilon}
\def\qed{{\hfill $\Box$ \bigskip}}
\def\NN{{\mathcal N}}
\def\AA{{\mathcal A}}
\def\MM{{\mathcal M}}
\def\BB{{\mathcal B}}
\def\CC{{\mathcal C}}
\def\LL{{\mathcal L}}
\def\DD{{\mathcal D}}
\def\FF{{\mathcal F}}
\def\EE{{\mathcal E}}
\def\QQ{{\mathcal Q}}
\def\WW{{\mathcal W}}
\def\RR{{\mathbb R}}
\def\R{{\mathbb R}}
\def\L{{\bf L}}
\def\K{{\bf K}}
\def\S{{\bf S}}
\def\A{{\bf A}}
\def\E{{\mathbb E}}
\def\F{{\bf F}}
\def\P{{\mathbb P}}
\def\N{{\mathbb N}}
\def\eps{\varepsilon}
\def\wh{\widehat}
\def\wt{\widetilde}
\def\pf{\noindent{\bf Proof.} }
\def\pff{\noindent{\bf Proof} }
\def\cp{\mathrm{Cap}}

\title{\Large \bf {Unavoidable collections of balls for censored stable processes}}

\author{{\bf Ante Mimica}
\quad {\bf Zoran Vondra\v{c}ek}
}

\date{}

\maketitle

\begin{abstract}
We study avoidability of collections of balls in bounded $C^{1,1}$ opens sets for censored $\alpha$-stable processes, $\alpha\in (1,2)$. The results are analog to the ones obtained for Brownian motion in S.~J.~Gardiner, M.~Ghergu, \emph{Champagne subregions of the unit ball with unavoidable bubbles},  Ann.~Acad.~Sci.~Fenn.~Math.  {\bf 35}  (2010)  321–-329. On the way we derive a Wiener-Aikawa-type criterion for minimal thinness with respect to the censored stable processes.
\end{abstract}

\noindent {\bf AMS 2010 Mathematics Subject Classification}: Primary 60J45; Secondary 31B15.

\noindent {\bf Keywords and phrases}: Censored stable process, minimal thinness

\section{Introduction}\label{s:intro}

Let $A$ be a Borel subset of the unit ball $B(0,1)\subset \R^d$, $d\ge 2$, not containing the origin. Then $A$ is said to be unavoidable if Brownian motion starting from the origin almost surely hits $A$ before hitting the boundary $\partial B$. More precisely, let $X=(X_t,\P_x)$ denote a standard Brownian motion in $\R^d$ and let $T_C=\inf\{t>0:\, X_t\in C\}$ be the hitting time of a Borel set $C\subset \R^d$. Then $A\subset B$ is \emph{unavoidable} if $\P_0(T_A <T_{\partial B(0,1)})=1$ and \emph{avoidable} if $\P_0(T_A<T_{\partial B(0,1)})<1$.

Problem of avoidability for sets $A$ that are unions of  balls has been studied recently.  Let $\{\overline{B}(x_n,r_n)\}_{n\ge 1}$ be a collection of pairwise disjoint closed balls contained in $B(0,1)$ satisfying  $|x_n|\to 1$ and $\sup\limits_{n\ge 1}\frac{r_n}{1-|x_n|}<1$. Define $A:=\cup_{n=1}^{\infty} \overline{B}(x_n,r_n)$. The domain $B(0,1)\setminus A$ is often called the champagne region and the balls are called bubbles.
Avoidability of balls in the unit disc in $\R^2$ was studied in \cite{Ake,OCS}, and in higher dimensions in \cite{Don}. Those results were extended in \cite{GG}.

The aim of this paper is to prove analogous results for a class of jump processes.
Note that if Brownian motion is replaced by the rotationally invariant $\alpha$-stable process, then any collection of balls in the unit ball is avoidable since the process jumps out of the unit ball with positive probability before hitting $A$.

The natural choice for the jump process replacing Brownian motion in avoidability problems in $B(0,1)$ is the censored $\alpha$-stable process with $\alpha\in (1,2)$. Roughly, this process is constructed from the symmetric $\alpha$-stable process by suppressing the jumps landing outside of the ball and continuing at the place where the suppressed jump has occurred. Such a process is transient (for $\alpha\in (1,2)$) and converges to the boundary at its lifetime. In particular, with probability one it cannot be killed while inside the state space -- if this were possible every subset of $B$ would be avoidable.
The censored stable processes were rigorously constructed and studied in \cite{BBC}. Fine properties of their potential theory in bounded $C^{1,1}$ open sets was further developed in \cite{CK} (for more detail see Section \ref{s:prelim}). These results will allow  us to replace the unit ball by a bounded $C^{1,1}$ open set.

Let $D\subset \R^d$, $d\ge 2$, be a bounded $C^{1,1}$ open set, denote by $\delta_D(x)$ the distance of $x\in D$ to the boundary $\partial D$,  and let $Y=(Y_t,\P_x)$ be the censored $\alpha$-stable process in $D$ with $\alpha\in (1,2)$. The lifetime of $Y$ will be denoted by $\zeta$. By \cite[Theorem 1.1]{BBC} it holds that $Y_{\zeta-}\in \partial D$.

The concept of avoidability for the process $Y$ and the set $D$ is defined analogously to the case of the Brownian motion and the unit ball. To be more precise, let $\{\overline{B}(x_n,r_n)$ be a collection of pairwise disjoint closed balls in $D$ such that $\delta_D(x_n)\to 0$ and $\sup\limits_{n\ge 1} \frac{r_n}{\delta_D(x_n)}<1/2$, and let $A=\cup_{n\ge 1}\overline{B}(x_n,r_n)$. The open set $D\setminus A$ will be called a champagne subregion of $D$, and as before the balls are called bubbles. Fix a point $x_0\in D$ such that $x_0\notin A$, and let $T_A=\inf\{t>0:\, Y_t\in A\}$. We will say that $A$ is avoidable if $\P_{x_0}(T_A<\zeta)<1$, and unavoidable otherwise. The first main result of this paper is the analog of \cite[Theorem 1]{GG}.

\begin{theorem}\label{t:theorem-1}
Let $D\setminus A$ be the champagne subregion of a bounded $C^{1,1}$ open subset $D$.
\begin{itemize}
\item[(a)] If $A$ is unavoidable, then
\begin{equation}\label{e:thm-1-1}
\sum_{n\ge 1}\frac{\delta_D(x_n)^{2\alpha-2}}{|x_n-z|^{d+\alpha-2}}\, r_n^{d-\alpha}=\infty
\quad \text{for }\sigma-\textrm{a.e.~}z\in \partial D\, .
\end{equation}
\item[(b)] Conversely, if \eqref{e:thm-1-1} and the separation condition
\begin{equation}\label{e:thm-1-2}
\inf_{j\neq k} \frac{|x_j-x_k|}{r_k^{1-\alpha/d} \delta_D(x_k)^{\alpha/d}} >0\
\end{equation}
hold, then $A$ is unavoidable.
\end{itemize}
\end{theorem}
Here $\sigma$ denotes the surface measure on $\partial D$.

The analog of \cite[Theorem 2]{GG} seems to make sense only in the unit ball $B=B(0,1)$. It concerns radii $r_n$ which are of the form $r_n=(1-|x_n|)\phi(|x_n|)$ where $\phi:[0,1)\to (0,1)$ is decreasing. For $a\in (0,1)$ let
$$
N_a(x)=\sharp [B(x, a(1-|x|))\cap \{x_n:\, n\in \N\}]
$$
be the number of centers in the ball $B(x, a(1-|x|)$, and let $M:[0,1)\to [1,\infty)$ be an increasing function satisfying
\begin{equation}\label{e:condition-M}
M(1-\tfrac{t}{2})\le c M(1-t)\quad \text{for all }t\in (0,1)\, .
\end{equation}
\begin{theorem}\label{t:theorem-2}
Let $\phi:[0,1)\to (0,1)$ be a decreasing function and $M:[0,1)\to [1,\infty)$ be a function satisfying \eqref{e:condition-M}. Let $B\setminus A$ be a champagne subregion of the ball $B$ such that $r_n=(1-|x_n|)\phi(|x_n|)$.
\begin{itemize}
\item[(a)] If $A$ is unavoidable and there are constants $a\in (0,1)$ and $b>0$ such that $N_a(x)\le bM(|x|)$ for all $x\in B$, then
\begin{equation}\label{e:thm-2-1}
\int_0^1 \frac{\phi(t)^{d-\alpha}M(t)}{1-t}\, dt=\infty \, .
\end{equation}
\item[(b)] Conversely, if \eqref{e:thm-2-1} holds together with the separation condition
\begin{equation}\label{e:thm-2-2}
\inf_{m\neq n} \frac{|x_m-x_n|}{\phi(|x_n|)^{1-\alpha/d} (1-|x_n|)} >0\, ,
\end{equation}
and there are constants $a\in (0,1)$ and $b\ge 0$ such that $N_a(x)\ge bM(|x|)$ for all $x\in B$, then $A$ is unavoidable.
\end{itemize}
\end{theorem}

By comparing the statements of Theorems \ref{t:theorem-1} and \ref{t:theorem-2} with \cite[Theorems 1 an 2.2]{GG} one sees that they are identical except for parameter $\alpha$ replacing 2. In proving Theorems \ref{t:theorem-1} and \ref{t:theorem-2} we will closely follow the ideas from \cite{GG}. In order to implement those ideas we had to develop certain potential-theoretic results for censored stable processes which are standard in case of Brownian motion (or classical potential theory). The first such result is quasi-additivity of capacity related to the censored stable process, see Proposition \ref{p:quasi-additivity}. Here we follow the exposition from \cite[Part II, Section 7]{AE} and only indicate the necessary changes mostly related to the construction of a comparable measure. The second necessary ingredient is a Wiener-Aikawa-type condition for minimal thinness of a set near the boundary point, see Propositions \ref{p:minthin-criterion-1} and \ref{p:aikawa-thinness}. The third result is a modification of a Aikawa-Borichev-type quasi-additivity of capacity from \cite[Theorem 3]{AB}, see Proposition \ref{p:ab}. With these results at hand we prove Theorems \ref{t:theorem-1} and \ref{t:theorem-2} in Sections \ref{s:thm1} and \ref{s:thm2}. Scaling plays a significant role in the proof of Theorem \ref{t:theorem-1}. We note that in case of censored stable process scaling is more  delicate than in case of Brownian motion due to the fact that both the process $Y$ and all related potential-theoretic properties are confined to the state space $D$.

We end this introduction with a few remarks on other works about avoidability of sets for jump processes. In a recent preprint \cite{HN13} the authors study the question of smallness of unavoidable sets in the context of balayage spaces. Their examples include Brownian motion, symmetric $\alpha$-stable processes in $\R^d$, $\alpha\in (0,2)$, and censored $\alpha$-stable processes in bounded $C^{1,1}$-open sets, $\alpha\in (1,2)$. In the latter case their results do not overlap with ours -- the goal in \cite{HN13} is to construct unavoidable collections of balls, or more generally unavoidable sets, having certain smallness properties, while Theorems \ref{t:theorem-1} and \ref{t:theorem-2} give sufficient and necessary conditions for a given collection of balls to be unavoidable. Another work that treats unavoidable collections of balls in $\R^d$ for certain class of isotropic L\'evy processes is \cite{MV} with results more in the spirit of the current paper.

Notation and conventions about constants:  The constants $C_G, C_M, C_H, C$ that we introduce in the next section, as well as the constant $C_1$ introduced in Section \ref{s:thm1}, stay fixed throughout the paper. The other constants, denoted by lowercase letters $c_1,c_2,\dots$, appear only locally in the paper and their numbering starts afresh in each section. Throughout the paper we use the notation $f(r) \asymp g(r)$ as  $r \to a$ to denote that $f(r)/g(r)$ stays between two positive constants as $r \to a$.

\section{Preliminaries about censored stable processes}\label{s:prelim}
Let $Y=(Y_t,\P_x)$ be a censored $\alpha$-stable process, $\alpha\in (1,2)$, in a bounded $C^{1,1}$ open set $D\subset \R^d$. Such process has been studied in \cite{BBC, CK, CS, Kim}. We will list several properties proved in those papers.

We first note that it follows from \cite{BBC} that $Y$ is a transient Hunt process with finite lifetime $\zeta$, and $Y_{\zeta-}\in \partial D$. The Dirichlet form $(\EE^D, \FF^D)$ of $Y$ is given by
\begin{equation}\label{e:df}
\EE^D(u,v)=\frac12 \AA(d,-\alpha)\int_D\int_D \frac{(u(x)-u(y))(v(x)-v(y))}{|x-y|^{d+\alpha}}\, dy\, dx\, ,\quad u,v\in C_c^{\infty}(D)\, ,
\end{equation}
where $\FF^D$ is the closure of $C_c^{\infty}(D)$ under $\EE^D_1=\EE^D+\langle \cdot, \cdot \rangle_{L^2(D)}$. Here $\AA(d,-\alpha)$ is an explicit, but unimportant constant.

The following Hardy's inequality is proved in \cite[Corollary 2.4]{CS}: There exists $c=c(D,\alpha)>0$ such that
\begin{equation}\label{e:hardy-inequality}
\EE^D(u,u)\ge c \int_D \frac{u(x)^2}{\delta_D(x)^{\alpha}}\, dx \, ,\quad u\in \FF^D\, .
\end{equation}

Let $G^D(x,y)$ denote the Green function of $Y$. The existence and sharp two-sided estimates for $G^D$ are proved in \cite[Theorem 1.1]{CK}: There exists $C_G=C_G(D,\alpha)\ge 1$ such that
\begin{eqnarray}\label{e:gfe}
\lefteqn{C_G^{-1}\left(1\wedge \frac{\delta_D(x)^{\alpha-1}}{|x-y|^{\alpha-1}}\right) \left(1\wedge \frac{\delta_D(y)^{\alpha-1}}{|x-y|^{\alpha-1}}\right)\, |x-y|^{\alpha-d}\le G^D(x,y)}\nonumber \\
&\le &C_G \left(1\wedge \frac{\delta_D(x)^{\alpha-1}}{|x-y|^{\alpha-1}}\right) \left(1\wedge \frac{\delta_D(y)^{\alpha-1}}{|x-y|^{\alpha-1}}\right)\, |x-y|^{\alpha-d}\, ,\qquad x,y\in D\, .
\end{eqnarray}
The constant $C_G(D,\alpha)$ can be chosen to be domain translation and dilation invariant.

Let $x_0\in D$ and define the Martin kernel based at $x_0$ by
$$
M^D(x,y)=\frac{G^D(x,y)}{G^D(x_0,y)}\, , \quad x,y\in D\, .
$$
It is proved in \cite[Theorem 1.2]{CK} that for each $Q\in \partial D$ there exists the limit $M^D(x,Q):=\lim_{y\to Q}M^D(x,y)$ and $M^D$ is jointly continuous on $D\times \partial D$. Further, there exists $C_M=C_M(x_0,D,\alpha)$ such that
\begin{eqnarray}\label{e:mke}
C_M^{-1} \frac{\delta_D(x)^{\alpha-1}}{|x-z|^{d+\alpha-2}}\le M^D(x,z)\le C_M \frac{\delta_D(x)^{\alpha-1}}{|x-z|^{d+\alpha-2}}\, ,\quad x\in D, z\in \partial D\, .
\end{eqnarray}
Consequently, the Martin boundary $\partial_M D$ of $D$ and the minimal Martin boundary $\partial_m D$ (with respect to $Y$) can be both identified with the Euclidean boundary $\partial D$.

Let $U$ be an open subset of $D$. A function $u:D\to [0,\infty)$ is harmonic in $U$ with respect to $Y$ if for every open set $V$ such that $V\subset \overline{V}\subset U$ it holds that
$$
u(x)=\E_x[u(Y_{\tau_V})]\quad \text{for all }x\in V\, .
$$
Here $\tau_V=\inf\{t>0:\, Y_t\notin V\}$ is the exit time of $Y$ from $V$. The function $u$ is regular harmonic in $U$ if
$$
u(x)=\E_x[u(Y_{\tau_U})]\quad \text{for all }x\in U\, .
$$
It is well known that regular harmonic functions are harmonic. Moreover, for $y\in D$, $x\mapsto G^D(x,y)$ is harmonic in $D\setminus \{y\}$, and for every $\epsilon >0$ regular harmonic in $D\setminus \overline{B}(y,\epsilon)$. Further, by \cite[Theorem 3.2]{BBC}, harmonic functions satisfy Harnack inequality. More precisely, we can deduce from that result that there exists a constant $C_H=C_H(d,\alpha)>0$ such that for every ball $B(x,r)\subset D$ and every nonnegative function $u$ on $D$ which is harmonic in $B(x,r)$,
\begin{equation}\label{e:harnack-inequality}
\sup_{y\in B(x,r/2)} u(y) \le C_H \inf_{y\in B(x,r/2)} u(y)\, .
\end{equation}

Let $\cp^D$ denote the capacity with respect to $Y$. It is proved in \cite[(3.10)]{CK} that there exists a constant $C=C(D,\alpha)\ge 1$ such that for every ball $B(x,r)\subset D$ satisfying $B(x,2r)\subset D$,
\begin{equation}\label{e:capacity-estimate}
C^{-1} r^{d-\alpha} \le \cp^D B(x,r) \le C r^{d-\alpha}\, .
\end{equation}

We look now at some scaling properties related to $Y$. For $r>0$, let $r^{-1}D:=\{x\in \R^d:\, rx\in D\}$. By \cite[Remark 2.3]{CK}, $\{r^{-1}Y_{r^{\alpha}t}, \P_x\}$ has the same distribution as the censored stable process in $r^{-1}D$ started at the point $r^{-1}x$ and
$$
G^{r^{-1}D}(x,y)=r^{d-\alpha}G^D(rx,ry)\, ,\quad x,y \in r^{-1}D\, .
$$
A simple computation using \eqref{e:df} gives that
$$
\EE^{r^{-1}D}(u,u)=r^{-d+\alpha} \EE^D(\wh{u},\wh{u})\, ,\quad u\in \FF^{r^{-1}D}\, ,
$$
where $\wh{u}(x)=u(r^{-1}x)$, $x\in D$. Therefore
\begin{eqnarray}\label{e:capacity-scaling}
\cp^{r^{-1}D}(A)&=&\inf\{\EE^{r^{-1}D}(u,u):\, u\in \FF^{r^{-1}D}, u\ge 1 \textrm{ a.e.~on }A\} \nonumber \\
&=&r^{-d+\alpha} \inf\{\EE^D(\wh{u},\wh{u}):\, \wh{u}\in \FF^D,\, \wh{u}\ge 1 \textrm{ a.e.~on }rA\} \nonumber \\
&=&r^{-d+\alpha} \cp^D(rA)\, .
\end{eqnarray}

Finally, let $\WW$ denote the family of all excessive functions with respect to $Y$. It is proved in \cite[Corollary 6.4]{HN13}, that $(D, \WW)$ is a balayage space in the sense of \cite{BH}. This will allow us to freely use results from \cite{BH}. Recall that for $u\in \WW$ and $B\subset D$, $R_u^B=\inf\{w\in \WW:\, w\ge v \text{ on }B\}$ is the reduced function of $u$ onto $B$, while its lower-semicontinuous regularization $\wh{R}_u^B\in \WW$ is called the balayage of $u$ onto $B$. If $T_B=\inf\{t>0:\, Y_t\in B\}$, then $\wh{R}_u^B(x)=\E_x u(X_{T_B})$ giving the probabilistic interpretation of the balayage, see \cite[VI.4]{BH}.

\section{Quasi-additivity of capacity}\label{s:quasi}
The goal of this section is to prove that $\cp^D$ is quasi-additive with respect to a Whitney decomposition of $D$.

Let $\{Q_j\}_{j\ge 1}$ be the Whitney decomposition of $D$. For each $Q_j$ let $Q_j^*$ denote the double of $Q_j$ and let $x_j$ denote the center of $Q_j$. Then $\{Q_j, Q_j^*\}$ is a quasi-disjoint decomposition of $D$, cf.~\cite[pp.~146-147]{AE}. A kernel $k:D\times D\to [0,+\infty]$ is said to satisfy the Harnack property with respect to $\{Q_j, Q_j^*\}$, cf.~\cite[p.~147]{AE}, if
$$
k(x,y)\asymp k(x',y)\, \textrm{ for all }x,x'\in Q_j \textrm{  and all }y\in D\setminus Q_j^*\, ,
$$
for all cubes $Q_j$ (with constants not depending on the cube). One way to get such kernels is as follows. Suppose that $u:D\to [0,\infty)$ is a function satisfying the scale invariant Harnack inequality of the form
$$
\sup_{Q_j}u\le c_1 \inf_{Q_j}u\, \quad \textrm{for all } Q_j\, ,
$$
where $c_1$ does not depend on $Q_j$. Typical $u$'s are the constant function $u\equiv 1$ and $u=g$ where $g(x):=G^D(x,x_0)\wedge 1$, and $x_0\in D$ fixed (see \eqref{e:gfe}). Define the kernel $k:D\times D\to [0,\infty]$ by

$$
k(x,y):=\frac{G^D(x,y)}{u(x)u(y)}\, ,\quad x,y\in D\, ,
$$
Note that $x\mapsto G^D(x,y)$ is regular harmonic in $Q_j$ for every $y\in D\setminus Q_j^*$. Hence by the scale invariant Harnack inequality \eqref{e:harnack-inequality} and the assumption on $u$, we see that $k$ satisfies the Harnack property with respect to $\{Q_j, Q_j^*\}$. For a measure $\lambda$ on $D$ let $\lambda_u(dy):=\lambda(dy)/u(y)$. Then
$$
K\lambda (x):=\int_D k(x,y)\, \lambda(dy)=\int_D \frac{G^D(x,y)}{u(x)u(y)}\, \lambda(dy)=\frac{1}{u(x)}\int_D G^D(x,y)\, \frac{\lambda(dy)}{u(y)}=\frac{1}{u(x)}\, G^D \lambda_u(dy)\, .
$$
We define the capacity with respect to the kernel $k$ as follows:
\begin{eqnarray*}
\cp_u(E):=\inf\{\|\lambda\|: K\lambda \ge 1 \textrm{ on }E\}\, .
\end{eqnarray*}
For a compact set $F\subset D$, consider  the balayage $\wh{R}_u^F$. Being a potential,  $\wh{R}_u^F=G^D \lambda_u^F$ for a measure $\lambda_u^F$ supported in $F$. Define the energy of $F$ (with respect to $u$) as
$$
\gamma_u(F):=\int_D \int_D G^D(x,y) \lambda_u^F(dx)\, \lambda_u^F(dy) =\int_D G^D \lambda_u^F(x)\, \lambda_u^F(dx)=\EE_D(G^D \lambda_u^F,G^D \lambda_u^F)\, .
$$
This definition of energy is in the usual way extended first to open, and then to Borel subsets of $D$. By using the dual definition of capacity
\begin{equation}\label{e:fuglede-equality}
\cp_u(F)=\sup\{\mu(F):\, \mu(D\setminus F)=0, K\mu \le 1\}\, ,
\end{equation}
for compact subsets $F\subset D$, see e.g.~\cite[Th\'eor\`eme 1.1]{Fug}, it is standard to show that
$$
\gamma_u(E)=\cp_u(E)\, ,\quad E\subset D\, .
$$
Note that in case $u\equiv 1$, $\gamma_1(E)=\cp_1(E)=\cp^D(E)$.

A Borel measure $\sigma_u$ (defined on Borel subsets if $D$) is comparable to the capacity $\cp_u$ with respect to $\{Q_j\}$ if there exists $c_2>0$ such that
\begin{eqnarray*}
& &\sigma_u(Q_j)\asymp \cp_u(Q_j)\, ,\quad \textrm{for all }Q_j\, ,\\
& &\sigma_u(E)\le c_2 \cp_u(E)\, ,\quad \textrm{for all Borel }E\, .
\end{eqnarray*}

Define
$$
\sigma_u(E):=\int_E u(x)^2\delta_D(x)^{-\alpha}\, dx\, ,\quad E\subset D\, .
$$
We claim that $\sigma_u$ is comparable with $C_u$. Note that on $Q_j$ we have $u\asymp u(x_j)$, hence $\wh{R}_u^{Q_j}\asymp u(x_j)\wh{R}_1^{Q_j}$, implying $G^D\lambda_u^{Q_j}\asymp u(x_j)G^D \lambda_1^{Q_j}$ (everywhere on $D$). Therefore,
\begin{eqnarray*}
\gamma_u(Q_j)&=&\int_{Q_j} G^D\lambda_u^{Q_j}(x)\, \lambda_u^{Q_j}(dx) \asymp u(x_j)\int_{Q_j} G^D\lambda_1^{Q_j}(x)\, \lambda_u^{Q_j}(dx)\\
&=&u(x_j)\int_{Q_j} G^D\lambda_u^{Q_j}(x)\, \lambda_1^{Q_j}(dx)\asymp u(x_j)^2 \int_{Q_j}G^D\lambda_1^{Q_j}(x)\lambda_1^{Q_j}(dx)\\
&=&u(x_j)^2 \cp^D(Q_j)\, .
\end{eqnarray*}
On the other hand,
\begin{eqnarray*}
\sigma_u(Q_j)&=&\int_{Q_j} u(x)^2\delta_D(x)^{-\alpha}\, dx \asymp u(x_j)^2 (\mathrm{diam}Q_j)^{-\alpha} \, |Q_j|^d \\
&\asymp &  u(x_j)^2 (\mathrm{diam}Q_j)^{-\alpha}\, (\mathrm{diam}Q_j)^d \asymp u(x_j)^2  \cp^D(Q_j)\, ,
\end{eqnarray*}
where the last asymptotic equality follows from \eqref{e:capacity-estimate}.
Thus, $\gamma_u(Q_j)\asymp \sigma_u(Q_j)$. Further, for any Borel $E\subset D$ and compact $F\subset E$, by using \eqref{e:hardy-inequality}, we have
\begin{eqnarray*}
\gamma_u(E)&\ge & \gamma_u(F)=\EE(G^D \lambda_u^F, G^D \lambda_u^F)\ge c_3 \int_D (G^D\lambda_u^F)(x)^2 \delta_D(x)^{-\alpha}\, dx\\
&\ge &c_3 \int_F (G^D\lambda_u^F)(x)^2 \delta_D(x)^{-\alpha}\, dx = c_3 \int_F u(x)^2 \delta_D(x)^{-\alpha}\, dx =c_3 \sigma_u(F)\, .
\end{eqnarray*}
This proves that $\gamma_u(E)\ge c_3  \sigma_u(E)$.

Now we can invoke \cite[Theorem 7.1.3]{AE} and conclude that $\gamma_u=\cp_u$ is quasi-additive with respect to $\{Q_j\}$.
\begin{proposition}\label{p:quasi-additivity}
The Green energy $\gamma_u$ is quasi-additive with respect to $\{Q_j\}$:
$$
\gamma_u(E)\asymp \sum_{j\ge 1} \gamma_u(E\cap Q_j)\, .
$$
\end{proposition}

\section{Minimal thinness}\label{s:thin}
In this section we prove a Wiener-Aikawa-type conditions for minimal thinness of a set near the boundary point.

Recall that $M^D(x,z)$ denotes the Martin kernel at $z$ (based at $x_0\in D$). The Martin boundary $\partial_M D$ and the minimal Martin boundary $\partial_m D$ of $D$ (with respect to $Y$) are identified with its Euclidean boundary $\partial D$. Recall that a set $E\subset D$ is said to be minimally thin at $z\in  \partial_m D$ if $\wh R^E_{M^D(\cdot, z)}\neq M^D(\cdot, z)$, cf.~\cite{Fol}. It is known, see e.g.~\cite{KW}, that every excessive function $u$ of $Y$ can be uniquely represented as
$$
u(x)=G^D \mu(x)+M^D \nu(x)=\int_D G^D(x,y)\, \mu(dy) + \int_{\partial D}M^D(x,z)\, \nu(dz)\, .
$$
The function $M^D\nu$ is the greatest harmonic minorant of $u$.

By following the proof of \cite[Theorem 9.2.6]{AG} and using \cite[Lemma 2.7]{Fol} instead of \cite[Lemma 9.2.2(c)]{AG} one obtains the next characterizations of minimal thinness.

\begin{proposition}\label{p:thinness}
Let $A\subset D$. The following are equivalent:

\noindent (a)
$A$ is minimally thin at $z\in \partial_m D$;

\noindent (b) There exists an excessive function $u=G^D\mu +M^D \nu$  such that
\begin{equation}\label{e:thinness}
\liminf_{x\to z, x\in A} \frac{u(x)}{M^D(x,z)}>\nu(\{z\})\, ;
\end{equation}

\noindent (c) There exists a potential $u=G^D\mu$ such that
\begin{equation}\label{e:thinness2}
\liminf_{x\to z, x\in A} \frac{u(x)}{M^D(x,z)}=+\infty \, .
\end{equation}
\end{proposition}
\pf We sketch the proof following the proof of \cite[Theorem 9.2.6]{AG}. Clearly, (c) implies (b). Assume that (b) holds. Then there exists a Martin topology neighborhood $W$ of $z$ and $a>\nu(\{z\})$ such that $u\ge a M^D(\cdot, z)$ on $A\cap W$.
If $\wh{R}^{A\cap W}_{M^D(\cdot, z)} =M^D(\cdot, z)$, then $u\ge \wh{R}^{A\cap W}_u\ge a M^D(\cdot,z)$ everywhere. Thus $u-aM^D(\cdot, z)$ is excessive, hence $u-aM^D(\cdot, z)=G^D\mu +M^D \wt{\nu}$ for a (non-negative) measure $\wt{\nu}$ on $\partial D$. On the other hand, $u-aM^D(\cdot, z)=G^D\mu +M^D \nu_{|\partial D\setminus \{z\}}+(\nu(\{z\})-a)M^D(\cdot, z)$.
This implies that $\wt{\nu}=\nu_{|\partial D\setminus \{z\}}+(\nu(\{z\})-a)\delta_z$ yielding $\wt{\nu}(z)=\nu(\{z\})-a<0$, which is a contradiction. Hence $\wh{R}^{A\cap W}_{M^D(\cdot, z)} \neq M^D(\cdot,z)$, i.e., $A$ is minimally thin at $z$. Thus (b) implies (a).

Suppose that (a) holds. By \cite[Lemma 2.7]{Fol}, there exists an open subset $U\subset \R^d$ such that $A\subset U$, and $U$ is minimally thin at $z$. By the analog of \cite[Theorem 9.2.5]{AG}, there is a decreasing sequence $(W_n)_{n\ge 1}$ of Martin topology open neighborhoods of $z$ shrinking to $z$ and such that $\wh{R}^{U\cap W_n}_{M^D(\cdot,z)}(x_0)\le 2^{-n}$. Let $u_1:=\sum_{n=1}^{\infty} \wh{R}^{U\cap W_n}_{M^D(\cdot,z)}$. Then $u_1$ is a sum of potentials, hence a potential itself since $u_1(x_0)<\infty$. Further, $\wh{R}^{U\cap W_n}_{M^D(\cdot,z)}={M^D(\cdot,z)}$ on the open set $U\cap W_n$. Therefore, $u_1(x)/M^D(x,z)\to \infty$ as $x\to z$, $x\in U$. Thus (c) holds.
\qed

Let $D$ be a $C^{1,1}$ open set. Fix $x_0\in D$ and let $M^D$ be the Martin kernel of $D$ based at $x_0$. The following proposition is an analog of \cite[Proposition V. 4.15]{BH}. A similar result is proved in \cite{KSV} in case of isotropic L\'evy processes satisfying certain conditions.

\begin{proposition}\label{p:minthin-criterion-1}
Let $E\subset D$ and let $z\in \partial D$. Define
$$
E_n=E\cap \{x\in D:\, 2^{-n-1}\le |x-z|<2^{-n}\}\, ,\quad n\ge 1\, .
$$
Then $E$ is minimally thin at $z$ if and only if $\sum_{n=1}^{\infty}R^{E_n}_{M^D(\cdot, z)}(x_0)<\infty$.
\end{proposition}
\pf Assume that $\sum_{n=1}^{\infty}R^{E_n}_{M^D(\cdot, z)}(x_0)<\infty$. Then there exists $n_0\in \N$ such that
$$
\sum_{n=n_0}^{\infty}R^{E_n}_{M^D(\cdot, z)}(x_0)<\frac12 M^D(x_0,z)\, .
$$
Let $B=B(z,2^{-n_0})$. Then $A:=B\cap E\subset \cup_{n=n_0}^{\infty} E_n$. Therefore,
$$
R^A_{M^D(\cdot, z)}(x_0)\le \sum_{n=n_0}^{\infty}R^{E_n}_{M^D(\cdot, z)}(x_0)<\frac12 M^D(x_0,z)\, ,
$$
implying $\wh{R}^A_{M^D(\cdot, z)}(x_0)<\frac12 M^D(x_0,z)$.
Hence, there exists an excessive function $u$ such that $u\ge M^D(\cdot, z)$ on $A$ and $u(x_0)<\frac12 M^D(x_0,z)$. By the Riesz decomposition,
$u=G^D\mu +M^D\nu$. Thus, $M^D \nu(x_0)=\int_{\partial_M D} M^D(x_0,z)\, \nu(dz)< \frac12 M^D(x_0,z)$ implying that $\nu(\{z\})<\frac12$.  Therefore,
$$
\liminf_{x\to z, x\in A} \frac{u(x)}{M^D(x,z)}\ge 1 >\frac12 >\nu(\{z\})\, .
$$
By Proposition \ref{p:thinness}, $A$ is minimally thin at $z$. Clearly, $E$ is also minimally thin at $z$.

Conversely, suppose that $E$ is minimally thin at $z$. By Proposition \ref{p:thinness}, there exists a potential $u$ such that $$
\liminf_{x\to z, x\in E} \frac{u(x)}{M^D(x,z)}=+\infty\, .
$$
Let $c_1=\max\{C_G, C_M\}$ where $C_G$ and $C_M$ are constants from \eqref{e:gfe} and \eqref{e:mke}, respectively. Without loss of generality, we may assume that
\begin{equation}\label{e:bound-at-0}
u(x_0)\le \frac{1}{2c_2}\, ,\quad \textrm{where }
c_2:=8c_1^4  \left(\frac{8}{3}\right)^{d}\, .
\end{equation}
There exists $n_1\in \N$ such that $u(x)>M^D(x,z)$ for all $x\in E\cap B(z,2^{-n_1})$. Thus, $E\subset \overline{B}(z,2^{-n_1})^c\cup
\{u>M^D(\cdot,z)\}$. For $n\ge n_1$ define
$$
G_n=\left\{x\in D:\, 2^{-n-1}<|x-z|<2^{-n}, u(x)>M^D(x,z)\right\}\qquad \textrm{and}\qquad G=\bigcup_{n=n_1}^{\infty}G_n\, .
$$
Let $x\in E_n$. Since $|x-z|\le 2^{-n_1}$, we have that $u(x)>M^D(x,z)$ and thus $x\in G_n$. This shows that $E_n\subset G_n$, $n\ge n_1$. Therefore, it suffices to show that $\sum_{n=n_1}^{\infty}R^{G_n}_{M^D(\cdot, z)}(x_0)<\infty$. Since $u>M^D(\cdot, z)$ on $G$,
it follows that $R^{G}_{M^D(\cdot, z)}(x_0)\le u(x_0)$.

Let $i\in\{1,2,3\}$. For every $n\in \N$ let
$$
U_n=G_{n_1+3n+i}\, .
$$
Since $i\in \{1,2,3\}$ is arbitrary, it suffices to show that $\sum_{n=1}^{\infty}R^{U_n}_{M^D(\cdot, z)}(x_0)<\infty$. Let $U=\bigcup_{n=1}^{\infty}U_n$. Then $U\subset G$ and thus $R^{U}_{M^D(\cdot, z)}(x_0)\le u(x_0)$. Note that since $U$ is open, $\wh{R}^U_{M^D(\cdot, z)}=R^U_{M^D(\cdot, z)}$ (see \cite[p.205]{BH}).
Since $u$ is a potential, the same holds for $\wh{R}^U_{M^D(\cdot, z)}$, hence there exists a measure $\mu$ such that $R^U_{M^D(\cdot, z)}=G^D\mu$. Moreover, since $R^U_{M^D(\cdot, z)}$ is harmonic on $\overline{U}^c$ (cf. \cite[III.2.5]{BH}), $\mu(U^c)=0$. Let $\mu_n:=\mu_{|\overline{U}_n}$. Since $\overline{U}_n$ are pairwise disjoint,
$$
\mu=\sum_{n=1}\mu_n \qquad \textrm{and} \qquad G^D\mu=\sum_{n=1}^{\infty}G^D \mu_n\, .
$$
Fix $n\in \N$ and consider $l\in \N$, $x\in U_n$, $y\in \overline{U}_l$. Then
$$
|x-y|\ge \frac34 |x-z|\, ,\qquad x\in U_n, y\in \overline{U_l}\, , l\neq n\, .
$$

Define $\mu_n'=\mu-\mu_n$ and let $x\in U_n$. By using \eqref{e:gfe} and \eqref{e:mke},
\begin{eqnarray}
G^D\mu_n'(x)&=&\int_D G^D(x,y)\, \mu_n'(dy)= M^D(x,z)\int_D \frac{G^D(x,y)}{M^D(x,z)}\, \mu_n'(dy) \label{e:key-step-1}\\
&\le & c_1^2 M^D(x,z) \int_D \frac{|x-y|^{\alpha-d}\frac{\delta_D(x)^{\alpha-1}}{|x-y|^{\alpha-1}}\frac{\delta_D(y)^{\alpha-1}}{|x-y|^{\alpha-1}}} {\frac{\delta_D(x)^{\alpha-1}}{|x-z|^{d+\alpha-2}}}\, \mu_n'(dy) \nonumber\\
&=& c_1^2 M^D(x,z) \int_D \left(\frac{|x-z|}{|x-y|}\right)^{d+\alpha-2} \delta_D(y)^{\alpha-1}\, \mu_n'(dy)\nonumber \\
&\le & c_1^2 \left(\frac43\right)^d M^D(x,z) \int_D \delta_D(y)^{\alpha-1}\, \mu_n'(dy)\, .\nonumber
\end{eqnarray}
We now compare $\delta_D(y)^{\alpha-1}$ with $G^D(x_0,y)$.
By choosing $n_1$ even larger, we can assume that $\frac12 \delta_D(x_0)\le |x_0-y|$, $\delta_D(y)\le |x_0-y|$ and $\frac12 |x_0-z|\le |x_0-y|\le 2 |x_0-z| $. Therefore,
$$
1\wedge\frac{\delta_D(y)^{\alpha-1}}{|x_0-y|^{\alpha-1}}=\frac{\delta_D(y)^{\alpha-1}}{|x_0-y|^{\alpha-1}} \quad \text{ and } \quad 1\wedge \frac{\delta_D(x_0)^{\alpha-1}}{|x_0-y|^{\alpha-1}}\ge 2^{1-\alpha} \frac{\delta_D(x_0)^{\alpha-1}}{|x_0-y|^{\alpha-1}} \, .
$$

Hence,
\begin{eqnarray}
G^D(x_0,y)&\ge & c_1^{-1} |x_0-y|^{\alpha-d} \left(1\wedge \frac{\delta_D(x_0)^{\alpha-1}}{|x_0-y|^{\alpha-1}}\right)\left(1\wedge \frac{\delta_D(y)^{\alpha-1}}{|x_0-y|^{\alpha-1}}\right) \nonumber \\
&\ge & c_1^{-1} |x_0-y|^{\alpha-d} 2^{1-\alpha}\frac{\delta_D(x_0)^{\alpha-1}}{|x_0-y|^{\alpha-1}}\frac{\delta_D(y)^{\alpha-1}}{|x_0-y|^{\alpha-1}} \nonumber \\
&=&c_1^{-1}  2^{1-\alpha} |x_0-y|^{-d-\alpha+2}\,  \delta_D(x_0)^{\alpha-1}\delta_D(y)^{\alpha-1} \nonumber \\
&\ge & c_1^{-1}  2^{1-\alpha} 2^{-d-\alpha+2} |x_0-z|^{-d-\alpha+2} \delta_D(x_0)^{\alpha-1}\delta_D(y)^{\alpha-1}\nonumber  \\
&\ge & c_1^{-2} 2^{3-2\alpha-d} M^D(x_0,z) \delta_D(y)^{\alpha-1} \ge  c_1^{-2} 2^{-1-d}\delta_D(y)^{\alpha-1}\, . \label{e:G-V-lower-bound}
\end{eqnarray}
Thus  we have that
\begin{eqnarray}
G^D\mu_n'(x)&\le& c_1^2 \left(\frac43\right)^d M^D(x,z) \int_D c_1^2 2^{d+3} G^D(x_0,y)\, \mu_n'(dy) \label{e:key-step-2}\\
&=& c_2 M^D(x,z) G^D\mu_n'(x_0)\le
c_2 M^D(x,z) G^D \mu(x_0)\nonumber \\
&\le &c_2 M^D(x,z) u(x_0) \le \frac12 M^D(x,z)\, .\nonumber
\end{eqnarray}
Since $G^D\mu_n+G^D\mu_n'=G^D\mu=R^U_{M^D(\cdot, z)}=M^D(\cdot, z)$ on $U$, it follows that $G^D\mu_n=M^D(\cdot, z)-G^D\mu_n'\ge M^D(\cdot, z)-\frac12 M^D(\cdot, z)=\frac12 M^D(\cdot, z)$ on $U_n$. This implies that $G^D\mu_n \ge \frac12 R^{U_n}_{M^D(\cdot, z)}$. Finally,
$$
\sum_{n=1}^{\infty}R^{U_n}_{M^D(\cdot, z)}(x_0)\le 2\sum_{n=1}^{\infty}G^D\mu_n(x_0)=2G^D\mu(x_0)<\infty\, .
$$
\qed

\medskip

Suppose that $\frac12 \delta_D(x_0)\le |x_0-y|\le 2\delta_D(x_0)$, $\delta_D(y)\le |x_0-y|$ and $|x_0-y|\ge \frac12 |x_0-z| $. It is shown in the proof of Proposition \ref{p:minthin-criterion-1}, see \eqref{e:G-V-lower-bound}, that $G^D(y,x_0)\ge c_3^{-1} \delta_D(y)^{\alpha-1}$ (with an explicit constant $c_3>1$). In the same way (even easier) it follows that $G^D(y,x_0)\le c_3 \delta_D(y)^{\alpha-1}$. Let $g(y):=G^D(y,x_0)\wedge 1$. It follows from the above discussion that $g(y)\asymp \delta_D(y)^{\alpha-1}$. Since
$$
M^D(y,z)\asymp \frac{\delta_D(y)^{\alpha-1}}{|y-z|^{d+\alpha-2}}\, ,
$$
we see that there exists $c_4>1$ such that
\begin{equation}\label{e:g-M-estimate}
c_4^{-1} \frac{g(x)}{|x-z|^{d+\alpha-2}} \le M^D(x,z)\le c_4 \frac{g(x)}{|x-z|^{d+\alpha-2}}\, ,\quad x\in D\, .
\end{equation}
In particular, if $E_n=E\cap \{x\in D:\, 2^{-n-1}\le |x-z|<2^{-n}\}$, then
$$
c_5^{-1}2^{n(d+\alpha-2)} g(x)\le M^D(x,z) \le c_5 2^{n(d+\alpha-2)} g(x)\, ,\quad x\in E_n\, .
$$
This implies that
$$
c_5^{-1}2^{n(d+\alpha-2)} R^{E_n}_g \le R^{E_n}_{M^D(\cdot, z)} \le c_5 2^{n(d+\alpha-2)} R^{E_n}_g\, .
$$
In particular,
\begin{equation}\label{e:equivalence-1}
\sum_{n=1}^{\infty}R^{E_n}_{M^D(\cdot, z)}(x_0)<\infty \quad \textrm{ if and only if}\quad \sum_{n=1}^{\infty}2^{n(d+\alpha-2)} R^{E_n}_g (x_0)<\infty\, .
\end{equation}

Note that $\wh{R}^{E_n}_g$ is a potential, hence there exists a measure $\lambda_n$ (supported by $\overline{E}_n$) such that $\wh{R}^{E_n}_g=G^D \lambda_n$. Also, $\wh{R}^{E_n}_g=g=G^D(\cdot, x_0)$ on $E_n$ (except a polar set, and at least for large $n$), hence
\begin{eqnarray*}
\wh{R}^{E_n}_g(x_0)&=& G^D \lambda_n(x_0)=\int_{\overline{E}_n}G^D(x_0,y)\, \lambda_n(dy)=\int_{\overline{E}_n}g(y)\, \lambda_n(dy)\\
&=&\int_{\overline{E}_n}\wh{R}^{E_n}_g(y)\, \lambda_n(dy)=\int_D \int_D G^D(x,y)\, \lambda_n(dy)\, \lambda_n(dx)= \gamma_g(E_n)\, ,
\end{eqnarray*}
the Green energy of $E_n$ with respect to $g$. We conclude from \eqref{e:equivalence-1} that
\begin{equation}\label{e:equivalence-2}
\sum_{n=1}^{\infty}R^{E_n}_{M^D(\cdot, z)}(x_0)<\infty \quad \textrm{ if and only if}\quad \sum_{n=1}^{\infty}2^{n(d+\alpha-2)} \gamma_g(E_n)<\infty\, .
\end{equation}
Thus we have proved the following Wiener-type criterion for minimal thinness.
\begin{corollary}\label{c:minthin-criterion-1}
Let $E\subset D$. Define
$$
E_n=E\cap \{x\in D:\, 2^{-n-1}\le |x-z|<2^{-n}\}\, ,\quad n\ge 1\, .
$$
Then $E$ is minimally thin at $z$ if and only if $\sum_{n=1}^{\infty}2^{n(d+\alpha-2)} \gamma_g(E_n)<\infty$.
\end{corollary}

Finally, we prove a version of Aikawa-type criterion for minimal thinness.
\begin{proposition}\label{p:aikawa-thinness}
Let $D$ be a bounded $C^{1,1}$ open set, $z\in \partial D$, $E\subset D$, and let $x_j$ denote the center of $Q_j$. The following are equivalent:

\begin{itemize}
\item[(a)] $E$ is minimally thin at $z$;
\item[(b)]
\begin{equation}\label{e:aikawa-thinness-2}
\sum_{j\ge 1} \mathrm{dist}(z,Q_j)^{-d-\alpha+2} g(x_j)^2 \cp^D(E\cap Q_j)<\infty\, ;
\end{equation}
\item[(c)]
\begin{equation}\label{e:aikawa-thinness}
\sum_{j\ge 1} \frac{\mathrm{dist}(Q_j,\partial D)^{2(\alpha-1)}}{\mathrm{dist}(z,Q_j)^{d+\alpha-2}} \, \cp^D(E\cap Q_j)<\infty\, .
\end{equation}
\end{itemize}
\end{proposition}
\pf (a)$\Leftrightarrow $ (b) By Corollary  \ref{c:minthin-criterion-1}, $E$ is minimally thin at $z$ if and only if $\sum_{n=1}^{\infty}2^{n(d+\alpha-2)} \gamma_g(E_n)<\infty$.  Further, let $A_n=\{x\in \R^d: 2^{-n-1}\le |x-z | <2^{-n}
\}$ so that $E_n=E\cap A_n$. If $A_n\cap Q_j \neq \emptyset$, then $\mathrm{dist}(z,Q_j)\asymp 2^{-n}$. By quasi-additivity of $\gamma_g$,
\begin{eqnarray*}
\sum_{n=1}^{\infty}2^{n(d+\alpha-2)} \gamma_g(E_n)&\asymp & \sum_{n=1}^{\infty}2^{-n(-d-\alpha+2)} \sum_{j\ge 1} \gamma_g(E_n\cap Q_j)\\
&\asymp &\sum_{j\ge 1} \sum_{n, A_n\cap Q_j\neq \emptyset} \mathrm{dist}(z,Q_j)^{-d-\alpha+2} \gamma_g(E_n\cap Q_j)\\
&=&\sum_{j\ge 1} \mathrm{dist}(z,Q_j)^{-d-\alpha+2}\sum_{n, A_n\cap Q_j\neq \emptyset}\gamma_g(E_n\cap Q_j)\\
&\asymp &\sum_{j\ge 1} \mathrm{dist}(z,Q_j)^{-d-\alpha+2} \gamma_g(E\cap Q_j)\, .
\end{eqnarray*}
For the last line we argue as follows: One inequality is subadditivity of capacity. For another note that there exists $N\in \N$ such that for every $Q_j$, $\sum_{n, A_n\cap Q_j\neq \emptyset}1=\sum_n 1_{A_n\cap Q_j}\le N$. Hence,
$\sum_{n, A_n\cap Q_j\neq \emptyset}\gamma_g(E\cap A_n\cap Q_j)\le \sum_{n, A_n\cap Q_j\neq \emptyset}\gamma_g(E \cap Q_j)\le N \gamma_g(E \cap Q_j)$.

Finally, we use that $\gamma_g(E\cap Q_j)\asymp g(x_j)^2 \cp^D(E\cap Q_j)$ which is proved in the same way as \cite[Part II, Proposition 7.3.1]{AE} to finish the proof.

\noindent
(b)$\Leftrightarrow $  (c) Note that
$$
g(x_j)=G^D(x_j, x_0)\wedge 1\asymp G^D(x_j, x_0)\asymp \delta_D(x_j)^{\alpha-1}\asymp \mathrm{dist}(Q_j,\partial D)^{\alpha-1}\, .
$$
hence the series  \eqref{e:aikawa-thinness-2} converges if and only if the series \eqref{e:aikawa-thinness} converges.
\qed

\section{Aikawa-Borichev-type result}\label{s:ab}
The purpose of this section is to prove the analog of \cite[Theorem 3]{AB} for the censored stable process.

Let $B(x,r)\subset D$ such that $B(x,2r)\subset D$. Define $\eta(r; x)>0$ as the radius such that
$$
|B(x, \eta(r;x))|=\cp^D B(x,r)\, ,
$$
and let $\eta^*(r; x)=\max\{\eta(r; x), 16r\}$. By \eqref{e:capacity-estimate} $C^{-1} r^{d-\alpha} \le \cp^D B(x,r) \le C r^{d-\alpha}$. Let $\sigma_d$ be the volume of the unit ball. Since $|B(x, \eta(r;x))|=\sigma_d \eta(r;x)^d$, we have that
$C^{-1} r^{d-\alpha} \le \sigma_d \eta(r;x)^d  \le C r^{d-\alpha}$, implying
\begin{equation}\label{e:estimate-eta}
C^{-1/d}\sigma_d^{-1/d}r^{1-\alpha/d}\le \eta(r;x) \le C^{1/d}\sigma_d^{-1/d}r^{1-\alpha/d}\, .
\end{equation}
Define
$$
\eta_l(r;x):=C^{-1/d}\sigma_d^{-1/d}r^{1-\alpha/d}\, \quad \textrm{and} \quad \eta_u(r;x):=C^{1/d}\sigma_d^{-1/d}r^{1-\alpha/d}\, .
$$

We recall now that a $C^{1,1}$ open set $D$ satisfies the following interior ball condition: There exists $R=R(D)>0$ such that for all $x_0\in D$ with $\delta_D(x_0)<R$ there is $z_0\in \partial D$ so that $|x_0-z_0|=\delta_D(x_0)$ and that $B(y_0,R)\subset D$ where $y_0=z_0+R(x_0-z_0)/|x_0-z_0|$. Note that if $a>0$, then $aD$ is again a $C^{1,1}$ open set and $R(aD)=aR(D)$.

\begin{lemma}\label{l:geometric}
Let $0<2r\leq r^*<R/2$ and let $x_0\in D$ be such that $\delta_D(x_0)<R/2$ and $B(x_0,r)\subset D$. There exists $\wt{x}\in D$ such that $|\wt{x}-x_0|=\frac34 r^*$ and for any $\theta\le \frac14$ it holds that
$$
B(\wt{x}, \theta r^*)\subset B(x_0,r^*)\cap D\, .
$$
Moreover, if $x\notin B(x_0,r^*)$ and $\rho=\mathrm{dist}(x, B(x_0,r))$, then $B(\wt{x}, \theta r^*)\subset B(x, 5\rho)$.
\end{lemma}
\pf
Let $z_0\in \partial D$ so that $|x_0-z_0|=\delta_D(x_0)$, and let $y_0=z_0+R(x_0-z_0)/|x_0-z_0|$. Then $B(y_0,R)\subset D$. Note that $r<\delta_D(x_0)$ since $B(x_0,r)\subset D$. If $y\in B(x_0,r)$, then
$$
|y-y_0|\le |y-x_0|+|x_0-y_0|<r +|z_0-y_0|-|z_0-x_0|=r+R-\delta_D(x_0)<R\, ,
$$
so $B(x_0,r)\subset B(y_0,R)$. Further,  $|x_0-y_0|=R-\delta_D(x_0)>R-R/2>r^*$. Let $\wt{x}$ be the point on the segment connecting $x_0$ and $y_0$ such that $|\wt{x}-x_0|=\frac{3}{4}r^*$.

Let $y\in B(\wt{x}, \theta r^*)$. Then
$$
|x_0-y|\le |x_0-\wt{x}|+\theta r^*\le \frac34 r^*+\frac14 r^*=r^*\, ,
$$
and
$$
|y_0-y|\le |y_0-\wt{x}|+|\wt{x}-y|=|y_0-x_0|-|x_0-\wt{x}|+|\wt{x}-y|<R-\delta_D(x_0)-\frac34 r^*+\theta r^*<R\, .
$$
This shows that $B(\wt{x}, \theta r^*)\subset B(x_0,r^*)\cap B(y_0,R)\subset  B(x_0,r^*)\cap D$.

Let $x\notin B(x_0,r^*)$; then $|x-x_0|= \rho+r \le \rho+r^*$ and $\rho\ge r^*-r>r^*/2$. If $y\in B(x_0,r^*)$, then
$$
|y-x|\le |y-x_0|+|x_0-x|<r^* +\rho + r^* \le 5\rho\, .
$$
Hence $B(x_0,r^*)\subset B(x, 5\rho)$, and thus also $B(\wt{x}, \theta r^*)\subset B(x, 5\rho)$.
\qed

\begin{proposition}\label{p:ab}
Suppose that $\{B(x_j,r_j)\}_{j\ge 1}$ is a collection of balls contained in $D$ such that $\delta_D(x_j)<R/2$ and $\eta^*(r_j; x_j)<R/2$ for all $j$, and the family $\{B(x_j, \eta^*(r_j; x_j))\}_{j\ge 1}$ is pairwise disjoint. If $E$ is a Borel set contained in $\cup_{j\ge 1} B(x_j, r_j)$, then
$$
\cp^D E\le \sum_{j\ge 1}\cp^D(B(x_j,r_j)\cap E)\le c \cp^D E\, ,
$$
where the constant $c>1$ depends on $Y$ only.
\end{proposition}

\pf
We follow the proof of \cite[Theorem 3]{AB} and indicate only necessary changes. For a finite measure $\nu$ on $D$, let $\|\nu\|=\nu(D)$ denote its total mass. First note that it follows from \eqref{e:fuglede-equality} that
$$
\cp^D (B(x_j,r_j)\cap E) =\sup\{\|\nu\|:\, \nu(D\setminus (B(x_j,r_j)\cap E))=0, G^D\nu \le 1 \text{ on }D\}\, .
$$
Hence, given $\epsilon >0$, there exists a measure $\mu_j$ concentrated on $B(x_j,r_j)\cap E$ such that $G^D\mu_j \le 1$ on $D$ and $\|\mu_j\|\ge \cp^D(B(x_j,r_j)\cap E)-2^{-j}\epsilon$. Let $\mu_j'$ be the measure on $D$ defined by
$$
d\mu_j'=\frac{\|\mu_j\|} {|B(x_j, \eta^*(r_j;x_j))|}\ind_{B(x_j, \eta^*(r_j;x_j))\cap D} \, dx\, .
$$
It is clear from the proof of \cite[Theorem 3]{AB} that once we prove the inequality
\begin{equation}\label{eq:ab_to_prove}
    G^D\mu_j(x)\leq  cG^D\mu_j'(x),\qquad x\in B(x_j,\eta^*(r_j;x_j))^c\, ,
\end{equation}
where $c=c(d,\alpha)>0$, the rest of the proof follows in the exactly same way.

For simplicity, we write $r_j^*:=\eta^*(r_j;x_j)$. We use Lemma \ref{l:geometric} with $x_j$ instead of $x_0$ and let $\wt{x}_j$ be the point corresponding to $\wt{x}$; then $|\wt{x}_j-x_j|=\frac{3}{4} r_j^*$ and with $\theta=1/8$ we have
$$
    B(\wt{x}_j,8^{-1}r_j^*)\subset B(x_j,r_j^*)\cap D\, ,
$$
Moreover, if  $x\in B(x_j,r_j^*)^c$ and $\rho_j=\mathrm{dist}(x, B(x_j,r_j))$, then $B(\wt{x}_j,8^{-1}r_j^*)\subset B(x,5\rho_j)$. In the following calculation we use that for $s>0$ and $t\ge 1$, it holds that
$$
1\wedge \frac{s}{t}\le 1\wedge s \le t\left(1\wedge \frac{s}{t}\right)\, .
$$
By \eqref{e:gfe} we have
 \begin{align}
 	G^D\mu_j'(x)&\geq C_G^{-1}\int_{B(x,5 \rho_j)\cap B(x_j,r_j^*)\cap D}\left(1\wedge \frac{\delta_D(x)^{\alpha-1}}{|x-y|^{\alpha-1}}\right)\left(1\wedge \frac{\delta_D(y)^{\alpha-1}}{|x-y|^{\alpha-1}}\right)|x-y|^{\alpha-d}\mu_j'(dy)\nonumber\\
	&\geq C_G^{-1}\left(1\wedge \frac{\delta_D(x)^{\alpha-1}}{5^{\alpha-1}\rho_j^{\alpha-1}}\right)5^{\alpha-d}\rho_j^{\alpha-d} \int_{B(\wt{x}_j,8^{-1}r_j^*)}\left(1\wedge \frac{\delta_D(y)^{\alpha-1}}{5^{\alpha-1}\rho_j^{\alpha-1}}\right)\mu_j'(dy)\nonumber\\
	&\geq C_G^{-1}5^{2-\alpha-d}\left(1\wedge \frac{\delta_D(x)^{\alpha-1}}{\rho_j^{\alpha-1}}\right)\rho_j^{\alpha-d} \inf_{y\in B(\wt{x}_j,8^{-1}r_j^*)}\left(1\wedge \frac{\delta_D(y)^{\alpha-1}}{\rho_j^{\alpha-1}}\right)\mu_j'(B(\wt{x}_j,8^{-1}r_j^*))\,.\label{eq:ab_tmp_1}
 \end{align}
 Similarly,
 \begin{align}
 	G^D\mu_j(x)&\leq C_G\int_{B(x_j,r_j)}\left(1\wedge \frac{\delta_D(x)^{\alpha-1}}{|x-y|^{\alpha-1}}\right)\left(1\wedge \frac{\delta_D(y)^{\alpha-1}}{|x-y|^{\alpha-1}}\right)|x-y|^{\alpha-d}\mu_j(dy)\nonumber\\
	&\leq C_G\left(1\wedge \frac{\delta_D(x)^{\alpha-1}}{\rho_j^{\alpha-1}}\right)\rho_j^{\alpha-d}\int_{B(x_j,r_j)}\left(1\wedge \frac{\delta_D(y)^{\alpha-1}}{\rho_j^{\alpha-1}}\right)\mu_j(dy)\nonumber\\
	&\leq C_G\left(1\wedge \frac{\delta_D(x)^{\alpha-1}}{\rho_j^{\alpha-1}}\right)\rho_j^{\alpha-d} \sup_{y\in B(x_j,r_j)}\left(1\wedge \frac{\delta_D(y)^{\alpha-1}}{\rho_j^{\alpha-1}}\right)\mu_j(B(x_j,r_j))\,.\label{eq:ab_tmp_2}
 \end{align}
  Let $y_1\in B(x_j,r_j)$ and $y_2\in B(\wt{x}_j,8^{-1}r_j^*)$. Then
 \[
 	|y_1-y_2|\leq r_j+\frac{3r_j^*}{4}+\frac{r_j^*}{8}\leq \frac{r_j^*}{16}+\frac{7r_j^*}{8}=\frac{15r_j^*}{16}\leq 15(r_j^*-r_j)\leq 15\delta_D(y_2)
 \]
 implying
 \[
 	\delta_D(y_1)\leq \delta_D(y_2)+|y_1-y_2|\leq 16\delta_D(y_2)\,.
 \]
 Hence, by (\ref{eq:ab_tmp_1}), (\ref{eq:ab_tmp_2}) and the last display,
 \begin{align*}
 	G^D\mu_j(x)&\leq C_G^2 5^{d+\alpha-2} 16^{\alpha-1} \frac{\mu_j(B(x_j,r_j))}{\mu_j'(B(\wt{x}_j,8^{-1}r_j^*))}G^D\mu_j'(x)\\
    &=C_G^2 5^{d+\alpha-2} 16^{\alpha-1}\frac{|B(x_j,r_j^*)|}{|B(\wt{x}_j,8^{-1}r_j^*)|}G^D\mu_j'(x)\\
	&=C_G^2 5^{d+\alpha-2} 16^{\alpha-1}8^d G^D \mu_j'(x)\,.
 \end{align*}
\qed

\section{Proof of Theorem \ref{t:theorem-1}}\label{s:thm1}
The constant function 1 is harmonic with respect to $Y$, hence there exists a measure $\mu$  on $\partial D$ such that
$$
1=\int_{\partial D}M^D(x,z)\, \mu(dz)\, .
$$
Since $M^D(x_0,z)=1$ for all $z\in \partial D$, we have that $1=\int_{\partial D}\mu(dz)=\|\mu\|$. Let $E\subset D$: Then
$$
R^E_1(x_0)=R^E_{\ \int_{\partial D}M^D(\cdot,z)\, \mu(dz)} (x_0)=\int_{\partial D} R^E_{M^D(\cdot, z)}(x_0)\, \mu(dz)\, .
$$
Since $R^E_{M^D(\cdot, z)}(x_0)\le M^D(x_0,z)=1$, we see that $R^E_1(x_0)=1$ if and only if $R^E_{M^D(\cdot, z)}(x_0)=1$ for $\mu$-a.e.~$z\in \partial D$. Thus, $R^E_1(x_0)=1$ if and only if $E$ is not minimally thin at $z$ for $\mu$-a.e.~$z\in \partial D$.

We note further that by \cite[Lemma 3.1]{Kim} $\mu$ is in fact the harmonic measure for $Y$ in $D$: $\mu(dz)=\P_{x_0}(Y_{\zeta-}\in dz)$. It is proved in \cite[Theorem 3.14]{Kim} that the harmonic measure $\mu$ is mutually absolutely continuous with respect to the surface measure $\sigma$ on $\partial D$. We conclude that $R^E_1(x_0)=1$ if and only if $E$ is not minimally thin at $z$ for $\sigma$-a.e.~$z\in \partial D$.

Let $\{\overline{B}(x_k,r_k)\}_{k\ge 1}$ be a family of disjoint closed balls in $D$, and let $A:=\cup_{k\ge 1}\overline{B}(x_k,r_k)$. Then  the family of  balls is unavoidable if $R^A_1(x_0)=1$, or, equivalently, if $A$ is not minimally thin at $z$ for $\sigma$-a.e.~$z\in \partial D$. As before, let $\{Q_j\}_{j\ge 1}$ be a Whitney decomposition of $D$. We will need the following simple geometric lemma whose proof is omitted.

\begin{lemma}\label{l:geometric-lemma}
Assume that $\sup_{k\ge 1} r_k /\delta_D(x_k)<1$. There exists a constant $C_1\ge 1$ such that for any $Q_j$ and any $\overline{B}(x_k,r_k)$ which intersects $Q_j$
\begin{equation}\label{e:geometric}
C_1^{-1}\le \frac{\mathrm{dist}(Q_j, \partial D)}{\delta_D(x_k)}\le C_1 \quad \textrm{and}\quad C_1^{-1}\le \frac{\mathrm{dist}(z,Q_j)}{|x_k-z|}\le C_1 \quad \textrm{for all }z\in \partial D\, .
\end{equation}
\end{lemma}

We first note that the number of cubes $Q_j$ which intersect a given ball $\overline{B}(x_k,r_k)$ is bounded above by a constant $c_2$ (independent of $k$).
Next note that if $B(x,r)$ is an open ball and $Q$ a closed cube, then $\cp^D(B(x,r)\cap Q)=\cp(\overline{B}(x,r)\cap Q)$. Indeed, every point in $\partial B(x,r)\cap Q$ is regular for $B(x,r)\cap Q$, hence $\P_y(T_{B(x,r)\cap Q}<\infty)=\P_y(T_{\overline{B}(x,r)\cap Q}<\infty)$
for all $y\in \R^d$. This shows that $B(x,r)\cap Q$ and $\overline{B}(x,r)\cap Q$ have the same capacitary measures, hence equal capacities.

Recall that $R=R(D)>0$ is the constant from the interior ball condition.

\begin{lemma}\label{l:aux}
Let $\{\overline{B}(x_k,r_k)\}_{k\ge 1}$ be a family of closed balls in $D$ such that $\sup r_k /\delta_D(x_k)<1/2$, $\delta_D(x_k)<R/2$ and $\eta^*(r_k; x_k)<R/2$ for all $k\ge 1$.   Assume that
\begin{itemize}
    \item[(i)] $r_k\le
    (16^dC \sigma_d)^{-1/\alpha}$ for all $k\ge 1$;
    \item[(ii)]
    $\displaystyle{\frac{|x_j-x_k|}{r_k^{1-\alpha/d}}\ge 2 C^{1/d}\sigma_d^{-1/d}}$, $ j\neq k$.
\end{itemize}
Then there exists $c_3=c_3(d)>0$ such that for every $j\ge 1$,
$$
\cp^D(A\cap Q_j)\ge c_3 \sum_{k\ge 1} \cp^D(\overline{B}(x_k,r_k)\cap Q_j)\, .
$$
\end{lemma}
\pf The first condition
$r_k\le (16^dC \sigma_d)^{-1/\alpha}$
is equivalent to
$16r_k\le C^{-1/d}\sigma_d^{-1/d}r_k^{1-\alpha/d}=\eta_l(r_k;x_k)\le \eta(r_k;x_k)$ which implies that $\eta^*(r_k; x_k)= \eta(r_k;x_k)$. The second condition implies that
$$
|x_j-x_k|\ge 2 C^{1/d}\sigma_d^{-1/d}r_k^{1-\alpha/d}=2 \eta_u(r_k;x_k)\ge 2\eta(r_k;x_k)\, ,
$$
that is, the balls $B(x_j, \eta^*(r_j;x_j))$ are disjoint. The claim now follows from Proposition \ref{p:ab} and the fact that $\cp^D(B(x_k,r_k)\cap Q_j)=\cp^D(\overline{B}(x_k,r_k)\cap Q_j)$. \qed

\begin{lemma}\label{l:aux-2}
Let $\{\overline{B}(x_k,r_k)\}_{k\ge 1}$ be a family of closed balls in $D$ such that $\sup r_k /\delta_D(x_k)<1/2$, $\delta_D(x_k)<R/2$ and $\eta^*(r_k; x_k)<R/2$ for all $k\ge 1$. Assume that
\begin{equation}\label{e:separation-2}
\frac{|x_j-x_k|}{r_k^{1-\alpha/d} \delta_D(x_k)^{\alpha/d}}\ge
32 C^{2/d}C_1^{2\alpha/d}\, ,\quad j\neq k\, ,
\end{equation}
where $C_1\geq 1$ is the constant from Lemma \ref{l:aux}. Then there exists a constant $c_4=c_4(d, D)$ such that
\begin{equation}\label{e:aux-2-aux}
\cp^D(A\cap Q_j)\ge c_4 \sum_{k\ge 1} \cp^D(\overline{B}(x_k,r_k)\cap Q_j)\, .
\end{equation}
\end{lemma}
\pf Let $Q_j$ be a Whitney cube such that
\begin{equation}\label{e:aux-2-aux-1}
\mathrm{dist}(Q_j, \partial D)\le C_1^{-1}(16^d C \sigma_d)^{-1/\alpha}\, .
\end{equation}
Define the scaling constant $a >0$ by
$$
a:=\frac{\left(16^d C\sigma_d\right)^{-1/\alpha}}{C_1 \mathrm{dist} (Q_j,\partial D)}\, .
$$
By \eqref{e:aux-2-aux-1} it holds that $a\ge 1$.

Let $y_k=a x_k$, $\rho_k=a r_k$ and $F:=\left(\cup_{k=1}^\infty B(y_k,\rho_k)\right)\cap (a Q_j)$. Since $\delta_{aD}(y_k)=a\delta_D(x_k)$, the condition $\sup \rho_k /\delta_{aD}(y_k)<1/2$ is satisfied. Let $R^a:=R(aD)=aR$. Then  $\delta_{aD}(y_k)=a\delta_D(x_k)<aR/2=R^a/2$.
We now check that $\eta^*_a(\rho_k;y_k)\le R^a/2$ where $\eta_a(\rho_k; y_k)$ is computed with respect to $\cp^{aD}$, the capacity with respect to the censored stable process in $aD$. Indeed, by \eqref{e:capacity-scaling}
$$
|B(y_k, \eta_a(\rho_k; y_k))|=\cp^{aD}B(y_k,\rho_k)=a^{d-\alpha}\cp^D B(x_k,r_k)=a^{d-\alpha}|B(x_k, \eta(r_k;x_k))|\, .
$$
This implies that $\eta_a(\rho_k;y_k)=a^{1-\alpha/d}\eta(r_k;x_k)\le a \eta(r_k;x_k)\le aR/2=R^a/2$ since $a\ge 1$. Clearly, $16\rho_k=16a r_k\le aR/2=R^a/2$.

Suppose that $B(x_k,r_k)\cap Q_j\neq \emptyset$. Then
$$
\rho_k=a r_k\le a \delta_D(x_k)\le a C_1 \mathrm{dist}(Q_j, \partial D)=
\left(16^d C \sigma_d\right)^{-1/\alpha}\, .
$$
Further, for $l\neq k$,
\begin{eqnarray*}
\frac{|y_l-y_k|}{\rho_k^{1-\alpha/d}}&=&\frac{a|x_l-x_k|}{a^{1-\alpha/d}\, r_k^{1-\alpha/d}}=a^{\alpha/d}\, \frac{|x_l-x_k|}{r_k^{1-\alpha/d}}\\
&\ge & 32 C^{2/d}C_1^{2\alpha/d}\delta_D(x_k)^{\alpha/d}a^{\alpha/d}\\
&=&2 C^{1/d} \sigma_d^{-1/d} \left(C_1 \frac{\delta_D(x_k)}{\mathrm{dist}(Q_j,\partial D)}\right)^{\alpha/d}\\
&\ge &2 C^{1/d} \sigma_d^{-1/d}\, .
\end{eqnarray*}
By Lemma \ref{l:aux}, there exists $c_3=c_3(d)$ such that
$$
\cp^{aD}(F\cap a Q_j)\ge c_3 \sum_{k\ge 1} \cp^{aD}(\overline{B}(y_k,\rho_k)\cap a Q_j)\, ,
$$
where $\cp^{aD}$ is the capacity with respect to the censored $\alpha$-stable process in $a D$. By \eqref{e:capacity-scaling}, $\cp^{aD}(F\cap a Q_j)=\cp^{aD}(a(A\cap Q_j))=a^{d-\alpha}\cp^D(A\cap Q_j)$ and $\cp^{aD}(\overline{B}(y_k,\rho_k)\cap a Q_j)=a^{d-\alpha}\cp^D(\overline{B}(x_k, r_k)\cap Q_j)$. Therefore
$$
\cp^D(A\cap Q_j)\ge c_3 \sum_{k\ge 1} \cp^D(\overline{B}(x_k,r_k)\cap Q_j)\, .
$$

For finitely many Whitney cubes that do not satisfy \eqref{e:aux-2-aux-1}, one obtains inequalities \eqref{e:aux-2-aux} by choosing $c_4\le c_3$ small enough.
\qed

\noindent
\emph{Proof of Theorem \ref{t:theorem-1}.}
(a) Assume that $A$ is unavoidable. Since $2r_k<\delta_D(x_k)$, we have that $B(x_k,2r_k)\subset D$. Hence, by \eqref{e:capacity-estimate}, $\cp^D \overline{B}(x_k,r_k)\le C r_k^{d-\alpha}$. If $\overline{B}(x_k,r_k)$ and $Q_j$ intersect, we have that
$$
\cp^D(\overline{B}(x_k,r_k)\cap Q_j)\le \cp^D(\overline{B}(x_k,r_k))\le C r_k^{d-\alpha}\, .
$$
Since the number of cubes $Q_j$ which intersect a given ball $\overline{B}(x_k,r_k)$ is bounded above by a constant $c_2$, we have for every $z\in \partial D$
\begin{eqnarray*}
\lefteqn{\sum_{j\ge 1} \frac{\mathrm{dist}(Q_j,\partial D)^{2(\alpha-1)}}{\mathrm{dist}(z,Q_j)^{d+\alpha-2}}\,  \cp^D(A\cap Q_j)}\\
&\le &\sum_{k\ge 1}\sum_{j\ge 1} \frac{\mathrm{dist}(Q_j,\partial D)^{2(\alpha-1)}}{\mathrm{dist}(z,Q_j)^{d+\alpha-2}}\, \cp^D(\overline{B}(x_k,r_k)\cap Q_j)\\
&\le & C_1^{2\alpha-2+d+\alpha-2} C c_2\sum_{k\ge 1}\frac{\delta_D(x_k)^{2\alpha-2}}{|x_k-z|^{d+\alpha-2}}\, r_k^{d-\alpha}\, .
\end{eqnarray*}
The claim now follows from \eqref{e:aikawa-thinness} in Proposition \ref{p:aikawa-thinness}\,.

(b) Conversely, assume that \eqref{e:thm-1-1} and the separation condition \eqref{e:thm-1-2} hold true. Consider only the balls $B(x_k,r_k)$ such that $\delta_D(x_k)<R/2$. In this way a finite number of balls is omitted. If we show that this smaller collection is unavoidable, the same will be true for the whole collection.

Choose $\delta\in (0,1]$ small enough so that
$$
\frac{|x_j-x_k|}{(\delta r_k)^{1-\alpha/d} \delta_D(x_k)^{\alpha/d}}\ge 32 C^{2/d}C_1^{2\alpha/d}\, ,\quad j\neq k\, ,
$$
and so that $\eta^*(\delta r_k; x_k)<R/2$.
Note that the latter is possible because
$$
\eta^*(\delta r_k;x_k)\le \max\left(C^{1/d}\sigma_d^{-1/d}\delta^{1-\alpha/d} r_k^{1-\alpha/d}, 16 \delta r_k\right)\, .
$$
 Let $A_\delta:=\cup_{k\ge 1}\overline{B}(x_k,\delta r_k)$.
Lemma \ref{l:aux-2} applied to the family of balls $\{B(x_k,\delta r_k)\}$  gives that
$$
\cp^D(A_\delta \cap Q_j) \ge c_4 \sum_{k\ge 1} \cp^D(\overline{B}(x_k,\delta r_k)\cap Q_j)\, .
$$
Combined with \eqref{e:geometric}  this yields
\begin{eqnarray*}
\lefteqn{\sum_{j\ge 1}\frac{\mathrm{dist}(Q_j, \partial D)^{2(\alpha-1)}}{\mathrm{dist}(z,Q_j)^{d+\alpha-2}}\, \cp^D(A_{\delta}\cap Q_j) }\\
&\ge &\sum_{j\ge 1}\sum_{k\ge 1}\frac{\mathrm{dist}(Q_j, \partial D)^{2(\alpha-1)}}{\mathrm{dist}(z,Q_j)^{d+\alpha-2}}\, \cp^D(\overline{B}(x_k,\delta r_k)\cap Q_j)\\
&\ge &\sum_{k\ge 1}\sum_{j\ge 1}\frac{c_4}{C_1^{d+3\alpha-4}}\, \frac{\delta_D(x_k)^{2\alpha-2}}{|x_k-z|^{d+\alpha-2}}\, \cp(\overline{B}(x_k,\delta r_k)\cap Q_j)\, .
\end{eqnarray*}
Subadditivity of $\cp^D$ implies that for each $k$ there is $j$ such that
$$
\cp^D(\overline{B}(x_k,\delta r_k)\cap Q_j)\ge c_2^{-1}\cp^D \overline{B}(x_k,\delta r_k)\ge c_2^{-1}C^{-1} \delta^{d-\alpha}r_k^{d-\alpha}\, ,
$$
where $c_2$ is the constant from the proof of part (a) and the last inequality follows from \eqref{e:capacity-estimate}.
Therefore,
\begin{eqnarray*}
\lefteqn{\sum_{j\ge 1}\frac{\mathrm{dist}(Q_j, \partial D)^{2(\alpha-1)}}{\mathrm{dist}(z,Q_j)^{d+\alpha-2}}\, \cp^D(A_{\delta}\cap Q_j)}\\
&\ge & \frac{c_4 \delta^{d-\alpha}}{C_1^{d-\alpha}c_2 C}\sum_{k\ge 1}\frac{\delta_D(x_k)^{2\alpha-2}}{|x_k-z|^{d+\alpha-2}}\, r_k^{d-\alpha} =+\infty
\end{eqnarray*}
for $\mu$-a.e.~$z\in \partial D$. It follows from Proposition \ref{p:aikawa-thinness}, \eqref{e:aikawa-thinness}, that $A_{\delta}$ is unavoidable, hence the same is true for $A$. \qed

\section{Proof of Theorem \ref{t:theorem-2}}\label{s:thm2}
Let $B=B(0,1)$, $M^B(x,z)$  the Martin kernel based at $x_0=0$, and $\mu$ be the measure on $\partial B$ such that
$$
1=\int_{\partial B}M^B(x,z)\, \mu(dz)\, .
$$
It holds by \eqref{e:mke} that
$$
C_M^{-1} \frac{(1-|x|)^{\alpha-1}}{|x-z|^{d+\alpha-2}}\le M^B(x,z)\le C_M \frac{(1-|x|)^{\alpha-1}}{|x-z|^{d+\alpha-2}}\, ,\quad x\in B, z\in \partial B\, .
$$
The last two displays imply
\begin{equation}\label{e:mke-consequence}
\frac{C_M^{-1}}{(1-|x|)^{\alpha-1}}\le \int_{\partial B} \frac{\mu(dz)}{|x-z|^{d+\alpha-2}} \le \frac{C_M}{(1-|x|)^{\alpha-1}}\, .
\end{equation}
With this inequality at hand, the proof is essentially the same as the proof of \cite[Theorem 2]{GG}. We provide the proof for readers' convenience. Recall that $a,b,c$ are constants from the statement of the theorem.

\noindent
\emph{Proof of Theorem 1.2.}
(a) Without loss of generality we may assume that $\phi(|x_k|)\le 1/2$ for all $k\ge 1$, so that $\sup_{k\ge 1}r_k/\delta_B(x_k)=\sup_{k\geq 1}\phi(|x_k|)\le 1/2$.
It follows from Theorem \ref{t:theorem-1} (a) that
$$
\sum_{k\ge 1}\frac{\delta_B(x_k)^{2\alpha-2}}{|x_k-z|^{d+\alpha-2}}\, r_k^{d-\alpha} =\sum_{k\ge 1} \frac{(1-|x_k|)^{d+\alpha-2}}{|x_k-z|^{d+\alpha-2}}\, \phi(|x_k|)^{d-\alpha}
=\infty \quad \mu-\textrm{a.e.~}z\in \partial B
$$
Let the center $x_k$ belong to the Whitney cube $Q_m$. Then
$$
1-|x_k|\le \mathrm{dist} (Q_m,\partial B)+\mathrm{diam} (Q_m)\, ,
$$
implying (recall that $\mathrm{diam} (Q_m)\le \mathrm{dist} (Q_m,\partial B) \le 4 \mathrm{diam} (Q_m)$)
\begin{equation}\label{e:urf-2}
1-|x_k|\le 5\, \mathrm{diam} (Q_m)\, ,\qquad 1-|x_k|\le 2(1-|x|)\, , \quad x\in Q_m\, .
\end{equation}
By Lemma \ref{l:geometric-lemma},
$$
|z-x_k|\ge \frac{1}{C_1}\, \mathrm{dist} (z,Q_m)\ge \frac{1}{2C_1}\left(\mathrm{dist} (z,Q_m)+\mathrm{diam} (Q_m)\right)\ge \frac{|z-x|}{2C_1}\, ,\quad x\in Q_m\, .
$$
Together with the fact that $\phi$ is decreasing, this gives the estimate
$$
\frac{(1-|x_k|)^{d+\alpha-2}}{|z-x_k|^{d+\alpha-2}}\, \phi(|x_k|)^{d-\alpha}\le \left(10 C_1 \mathrm{diam} (Q_m)\right)^{d+\alpha-2}\, \frac{\phi((2|x|-1)^+)^{d-\alpha}}{|z-x|^{d+\alpha-2}}\, , \quad x\in Q_m\, .
$$
Note that the number of centers $x_k$ that belong to $Q_m$ is bounded from above by $c_1(a,c,d)b M((2|x|-1)^+)$ for every $x\in Q_m$. By using that $\sup_{x\in Q_m}\mathrm{diam} (Q_m)/ (1-|x|)<1$, it follows that
\begin{eqnarray*}
\lefteqn{\sum_{x_k\in Q_m} \frac{(1-|x_k|)^{d+\alpha-2}}{|x_k-z|^{d+\alpha-2}}\, \phi(|x_k|)^{d-\alpha}}\\
&\le & c_2(a,c,d)b C_1^{d+\alpha-2}(\mathrm{diam}(Q_m))^{\alpha-2} \int_{Q_m} \frac{\phi((2|x|-1)^+)^{d-\alpha}M((2|x|-1)^+)}{|z-x|^{d+\alpha-2}}\, dx\\
&\le & c_3(a,b,c,C_1,d) \int_{Q_m} \frac{\phi((2|x|-1)^+)^{d-\alpha}M((2|x|-1)^+)(1-|x|)^{\alpha-2}}{|z-x|^{d+\alpha-2}}\, dx
\end{eqnarray*}
By summing up over all Whitney cubes, we get,
$$
\sum_{k\ge 1} \frac{(1-|x_k|)^{d+\alpha-2}}{|x_k-z|^{d+\alpha-2}}\, \phi(|x_k|)^{d-\alpha}\le c_3(a,b,c,C_1,d) \int_B \frac{\phi((2|x|-1)^+)^{d-\alpha}M((2|x|-1)^+)(1-|x|)^{\alpha-2}}{|z-x|^{d+\alpha-2}}\, dx
$$
Since the left-hand side is infinite for $\mu$-a.e.~$z\in \partial B$, the same is true for the right-hand side. By integrating the right-hand side over $\partial B$ with respect to $\mu$, and by using \eqref{e:mke-consequence}, we get that
$$
\int_B \frac{\phi((2|x|-1)^+)^{d-\alpha}M((2|x|-1)^+)}{1-|x|}\, dx =\infty\, .
$$
By switching to polar coordinates (and using that $t^{d-1}$ is bounded near 1) this yields
$$
\int_{1/2}^1 \frac{\phi(2t-1)^{d-\alpha}M(2t-1)}{1-t}\, dt=\infty\, .
$$
The change of variables implies \eqref{e:thm-2-1}.

\noindent
(b)
Let $(t_i)_{i\ge 1}$ be a sequence defined by $t_i:=1-\frac{1}{2}\left(\frac{1-a}{1+a}\right)^i$\,.
Then $(t_i)_{i\ge 1}$ is increasing and $t_i<1$.
Let $z\in \partial B$ and define $z_i:=t_i z$. Then the balls $B(z_i, a(1-|z_i|))$ are pairwise disjoint. Indeed, this follows from the inequality $t_i+a(1-t_i)< t_{i+1}-a(1-t_{i+1})$. Further, for $x\in B(z_i, a(1-|z_i|))$ it holds that
$$
1-|x|\ge 1-t_i-a(1-t_i)=(1-a)(1-t_i)\, \quad \textrm{and} \quad |z-x|\le 1-t_i+a(1-t_i)=(1+a)(1-t_i)\, .
$$
Therefore,
\begin{eqnarray*}
\lefteqn{\sum_{k\ge 1} \frac{(1-|x_k|)^{2\alpha-2}}{|z-x_k|^{d+\alpha-2}}\, r_k^{d-\alpha} =\sum_{k\ge 1} \frac{(1-|x_k|)^{d+\alpha-2}}{|z-x_k|^{d+\alpha-2}}\, \phi(|x_k|)^{d-\alpha}}\\
&\ge & \sum_{i\ge 1} \sum_{\{k:\, x_k\in B(z_i, a(1-|z_i|))\}} \frac{(1-|x_k|)^{d+\alpha-2}}{|z-x_k|^{d+\alpha-2}}\, \phi(|x_k|)^{d-\alpha}\\
&\ge &\left(\frac{1-a}{1+a}\right)^{d+\alpha-2} \sum_{i\ge 1}b M(|z_i|)
\phi(1-(1-a)(1-t_i))^{d-\alpha}\\
&\ge &c_4(a,b,c,\gamma, d) \sum_{i\ge 1}
\phi(1-(1-a)(1-t_i))^{d-\alpha}M(1-(1-a)(1-t_{i+1}))\, ,
\end{eqnarray*}
where in the last line we have used
$
	|z_i|=t_i=1-\tfrac{1+a}{(1-a)^2}(1-a)(1-t_{i+1}),
$
\eqref{e:condition-M} and the  fact that $M$ increases to conclude that $M(|z_i|)\geq c^{-\lceil\log_2\frac{1+a}{(1-a)^2}\rceil}M(1-(1-a)(1-t_{i+1}))$\,.

Let $s_i:=1-(1-a)(1-t_i)$. Then $s_{i+1}-s_i=(1-a)(t_{i+1}-t_i)=\frac{2a}{1+a}(1-s_i)$.
Rewriting the last expression in terms of $s_i$ we get
\begin{eqnarray*}
\lefteqn{\sum_{k\ge 1} \frac{(1-|x_k|)^{2\alpha-2}}{|z-x_k|^{d+\alpha-2}}\, r_k^{d-\alpha}\ge c_4(a,b,c,\gamma, d) \sum_{i\ge 1} \phi(s_i)^{d-\alpha}M(s_{i+1}) }\\
&=&c_4(a,b,c,\gamma, d)\tfrac{1+a}{2a}\sum_{i=1}^{\infty}  \frac{\phi(s_i)^{d-\alpha}M(s_{i+1})}{1-s_i}\, (s_{i+1}-s_i)\\
&=&c_4(a,b,c,\gamma, d)\tfrac{1+a}{2a}\sum_{i=1}^{\infty}\int_{s_i}^{s_{i+1}} \frac{\phi(s_i)^{d-\alpha}M(s_{i+1})}{1-s_i}\, ds\\
&\ge &c_4(a,b,c,\gamma, d)\tfrac{1+a}{2a}\sum_{i=1}^{\infty}\int_{s_i}^{s_{i+1}} \frac{\phi(s)^{d-\alpha}M(s)}{1-s}\, ds\\
&=&c_4(a,b,c,\gamma, d)\tfrac{1+a}{2a}\int_{s_1}^1 \frac{\phi(s)^{d-\alpha}M(s)}{1-s}\, ds\, .
\end{eqnarray*}
We conclude from the assumption \eqref{e:thm-2-1} that the last display is equal to $+\infty$. Since the separation condition \eqref{e:thm-2-2} is the same as \eqref{e:thm-1-2}, the claim follows from Theorem \ref{t:theorem-1} (b). \qed

\begin{corollary}\label{c:regularly-spaced}
Let $\phi$ be as above and assume that the family of balls  $\{\overline{B}(x_k,r_k)\}_{k\ge 1}$ satisfies the separation condition
\begin{equation}\label{e:separation-strong}
\inf_{j\neq k} \frac{|x_j-x_k|}{1-|x_k|}>0\, ,
\end{equation}
and that for some $a\in (0,1)$, $N_a(x)\ge 1$ for every $x\in D$. Then the family  $\{\overline{B}(x_k,r_k)\}_{k\ge 1}$ is unavoidable if and only if
\begin{equation}\label{e:urf-2}
\int_0^1 \frac{\phi(t)^{d-\alpha}}{1-t}\, dt=\infty\, .
\end{equation}
\end{corollary}
\pf First note that $N_a(x)\ge 1$ for every $x\in D$ and the separation condition \eqref{e:separation-strong} imply that there exists $b\ge 1$ such that $1\le N_a(x)\le b$. Hence, by taking $M(t)=1$ for all $t\in [0,1)$, we see that $M(|x|)\le N_a(x)\le b M(|x|)$ for all $x\in D$. Moreover, since $0<\phi(|x|)^{1-\alpha/d}<1$, \eqref{e:separation-strong} implies the weaker separation condition \eqref{e:separation-2}. The statement now follows immediately from Theorem \ref{t:theorem-2}
\qed

\medskip
{\bf Acknowledgement:} We thank P.~Kim and R.~Song for allowing us to use some ideas from  \cite{KSV}.

\vspace{.1in}
\begin{singlespace}

\end{singlespace}

\end{doublespace}

\end{document}